\documentclass[11pt, reqno]{amsart}
\usepackage{amsmath,amsthm,amssymb,graphicx}
\theoremstyle{plain}
\newtheorem{thm}{Theorem}[section]
\newtheorem{prop}[thm]{Proposition}
\newtheorem{lem}[thm]{Lemma}
\newtheorem{cor}[thm]{Corollary}

\theoremstyle{definition}
\newtheorem{rem}[thm]{Remark}
\newtheorem{defn}[thm]{Definition}
\newtheorem{eg}[thm]{Example}
\numberwithin{equation}{section}

\newcommand{\bthm}{\begin{thm}}
\newcommand{\ethm}{\end{thm}}
\newcommand{\bprop}{\begin{prop}}
\newcommand{\eprop}{\end{prop}}
\newcommand{\bcor}{\begin{cor}}

\newcommand{\ecor}{\end{cor}}
\newcommand{\blem}{\begin{lem}}
\newcommand{\elem}{\end{lem}}
\newcommand{\bproof}{\begin{proof}}
\newcommand{\eproof}{\end{proof}}
\newcommand{\bca}{\begin{cases}}
\newcommand{\eca}{\end{cases}}
\newcommand{\brem}{\begin{rem}}
\newcommand{\erem}{\end{rem}}
\newcommand{\bpm}{\begin{pmatrix}}
\newcommand{\epm}{\end{pmatrix}}
\newcommand{\bdefn}{\begin{defn}}
\newcommand{\edefn}{\end{defn}}
\newcommand{\bsub}{\begin{subtitle}}
\newcommand{\esub}{\end{subtitle}}
\newcommand{\ben}{\begin{enumerate}}
\newcommand{\een}{\end{enumerate}}
\newcommand{\beg}{\begin{eg}}
\newcommand{\eeg}{\end{eg}}

\newcommand{\C}{\mathbb{C}}

\newcommand{\R}{\mathbb{R}}

\newcommand{\ca}{\mathcal{A}}
\newcommand{\ci}{\mathcal{I}}

\newcommand{\cg}{\mathcal{G}}

\newcommand{\cu}{\mathcal{U}}

\newcommand{\ck}{\mathcal{K}}

\newcommand{\ad}{\mathrm{ad}}

\def\ms{\medskip}
\def\bs{\bigskip}

\def\p{\partial}

\def\Im{{\rm Im\/}}

\def\I{{\rm I\/}}

\def\Ad{{\rm Ad\/}}
\def\ti{\tilde}

\def \a {\alpha}
\def \b {\beta}

\def \e {\epsilon}

\def \l {\lambda}

\def \n {\,\vert\,}

\def \o {\theta}
\def\w{\omega}
\def\W{\Omega}

\def\R{\mathbb{R} }
\def\C{\mathbb{C}}

\def\cg{{\mathcal {G}}}

\def\ck{{\mathcal {K}}}
\def\cl{{\mathcal{L}}}

\def\co{{\mathcal {O}}}

\def\cu{{\mathcal {U}}}

\def \li{\langle}
\def \ri{\rangle}

\begin{document}

\title[Curved flats and conservation laws]
{Curved flats, exterior differential systems,\\ and conservation laws}
\author{Chuu-Lian Terng$^*$}\thanks{$^*$Research supported
in  part by NSF Grant DMS-0306446}
\address{Northeastern University and UC Irvine}
\email{terng@neu.edu}
\author{ ErXiao Wang$^\dag$}\thanks{$^\dag$Research supported
in  part by Postdoctoral fellowship of MSRI}
\address{MSRI, Berkeley, CA 94720}
\email{ewang@mrsi.org}

\ms \hskip 3in \today

\begin{abstract} 

Let $\sigma$ be an involution of a real semi-simple Lie group $U$, $U_0$ the subgroup fixed by $\sigma$, and $U/U_0$ the corresponding symmetric space.  Ferus and Pedit called a submanifold $M$ of a rank $r$ symmetric space $U/U_0$ a {\it curved flat\/} if $T_pM$ is tangent to an $r$-dimensional flat of $U/U_0$ at $p$ for each $p\in M$. They noted that the equation for curved flats is an integrable system.  Bryant used the involution $\sigma$ to construct an involutive exterior differential system $\ci_\sigma$ such that integral submanifolds of $\ci_\sigma$ are curved flats.  Terng used $r$ first flows in the $U/U_0$-hierarchy of commuting soliton equations to construct the $U/U_0$-system.  She showed that the $U/U_0$-system and the curved flat system are gauge equivalent, used the inverse scattering theory to solve the Cauchy problem globally with smooth rapidly decaying initial data, used loop group factorization to construct infinitely many families of explicit solutions, and noted that many these systems occur as the Gauss-Codazzi equations for submanifolds in space forms.  The main goals of this paper are: (i) give a review of  these known results, (ii) use techniques from soliton theory to construct infinitely many integral submanifolds and conservation laws for the exterior differential system $\ci_\sigma$.   
\end{abstract}

\maketitle

\section{Introduction}

Let $G$ be a complex semi-simple Lie group, $\tau$ an involution of $G$ such that its differential at the identity $e$ is complex conjugate linear, and $\sigma$ an inovlution of $G$ such that the differential is complex linear.  Assume that 
\begin{align}\label{aa}
\tau\sigma=\sigma\tau.
\end{align}
Let $U$ be the fixed point set of $\tau$, i.e., a real form of $G$. We will still use $\tau, \sigma$ to denote $d\tau_e$ and $d\sigma_e$ respectively.   Let 
 $\cg, \cu$ denote the Lie algebras of $G$ and $U$ respectively. Since $\sigma$ and $\tau$ commute, $\sigma(U)\subset U$. So $\sigma| U$ is an involution of $U$.   Let $\cu_0, \cu_1$ denote the $+1, -1$ eigenspaces of $\sigma$ on $\cu$.  Then 
 $$[\cu_0,\cu_0]\subset \cu_0, \quad [\cu_0, \cu_1]\subset\cu_1, \quad [\cu_1, \cu_1]\subset \cu_0.$$
 The quotient space $U/U_0$ is a symmetric space, and the eigen-decomposition $\cu=\cu_0+\cu_1$ is called a {\it Cartan decomposition\/}.  
 
 Ferus and Pedit (\cite{FerPed96a}) called a submanifold $M$ of a rank $r$ symmetric space $U/U_0$ a {\it curved flat\/} if $T_pM$ is tangent to an $r$-dimensional flat of $U/U_0$ at $p$ for each $p\in M$. They noted that the equation for curved flats is an integrable system.  Bryant (\cite{Bry03}) used the involution $\sigma$ to construct a natural involutive exterior differential system $\ci_\sigma$ such that integral submanifolds of $\ci_\sigma$ in $U$ project down to curved flats in $U/U_0$.  Terng (\cite{Ter97}) used $r$ first flows in the $U/U_0$-hierarchy of commuting soliton equations to construct the $U/U_0$-system.  She showed that the $U/U_0$-system and the curved flat system are gauge equivalent, used the inverse scattering theory to solve the Cauchy problem globally with smooth rapidly decaying initial data (\cite{Ter97}), used loop group factorization to construct infinitely many families of explicit solutions (\cite{TerUhl00a}), and noted that many these systems occur as the Gauss-Codazzi equations for submanifolds in space forms (\cite{Ter97, Ter03}).    The main goals of this paper are: (i) review some of these known results, (ii) use techniques from soliton theory to construct infinitely many integral submanifolds and conservation laws for the exterior differential system $\ci_\sigma$.
 
 We review the definitions of these systems next. 
An element $a\in \cu_1$ is called {\it regular\/} if 
\ben
\item[(i)] $\ca:=\{y\in \cu_1\n [a,y]=0\}$ is a maximal abelian subspace in $\cu_1$, 
\item[(ii)] Ad$(U_0)(\ca)$ is open in $\cu_1$.  
\een 

Let $(\ , )$ be an ad-invariant, non-degenerate bilinear form on $\cu$.  Given a linear subspace $V$ of $\cu$ let
$$V^\perp=\{y\in \cu\n (y, V)=0\}.$$
Assume that $U/U_0$ has rank $r$.  Let  $\ca$ be a maximal abelian subspace in $\cu_1$, and let $a_1, \ldots, a_r\in \ca$ be regular and form  a basis of $\ca$.  
{\it The $U/U_0$-system} (cf. \cite{Ter97}) is the following PDE for $v:\R^r\to \cu_1\cap\ca^\perp$:
 \begin{equation}\label{ai}
 [a_i, v_{x_j}]-[a_j, v_{x_i}]= [[a_i, v], [a_j, v]], \quad 1\leq i\not=j\leq r.
 \end{equation} 
 These systems occur naturally in submanifold geometry.  For example, the Gauss-Codazzi equations for isometric immersions of space forms in space forms (\cite{FerPed96b, Ter97, BrDuPaTe02}), for isothermic surfaces in $\R^n$ (\cite{Bur00, BrDuPaTe02}), and for flat Lagrangian submanifolds in $\C^n$ or in  $\C P^n$ (\cite{Ter03}).  
 
 The $U/U_0$-system also arises naturally from soliton theory (cf. \cite{Ter97}).  In fact, 
given $1\leq i\leq r$, $b\in \ca$, and  a positive integer $j$, the  $(b,j)$-th soliton flow in the $U/U_0$-hierarchy  is a certain partial differential equation for $v:\R^2\to \cu_1\cap\ca^\perp$:
$$v_t= P_{b, j}(v).$$
For example, the second flow in the $SU(2)$-hierarchy is the NLS (non-linear Schr\"odinger equation), the third flow in the $SU(2)/SO(2)$-hierarchy is the modified KdV equation, and the first flow in the $SU(3)/SO(3)$-hierarchy is the $3$-wave equation.  
 The $U/U_0$-system \eqref{ai} is given by the collection of the $(a_j, 1)$- flows with $1\leq j\leq r$ in the $U/U_0$-hierarchy.

The {\it curved flat system associated to $U/U_0$} (cf. \cite{FerPed96b}) is the following first order system for $(A_1, \ldots, A_r):\R^r\to \prod_{i=1}^r \cu_1$:
 \begin{equation}\label{bd}
 \bca (A_i)_{x_j}= (A_j)_{x_i}, \quad i\not= j, \\
 [A_i, A_j]=0, \quad i\not=j.
 \eca
 \end{equation}
 It is known that solutions of the curved flat system give rise to curved flats in $U/U_0$ (\cite{FerPed96a}).
 
Let $\a=g^{-1}dg$ be the Maurer-Cartan form on $U$. Write 
$\a= \a_0+\a_1$ with respect to the Cartan decomposition $\cu=\cu_0+\cu_1$.  Let $\ci_\sigma$ be the exterior differential ideal generated by $\a_0$.  
It was observed by Bryant (\cite{Bry03}) that $(U,\ci_\sigma)$ is involutive and the PDE for the exterior differential system $\ci_\sigma$ is  the curved flat system \eqref{bd}.  
If $f:\co\to U$ is a maximal integral submanifold of  the exterior differential system (EDS) $(U,\ci_\sigma)$, then $f^*(\a_0)=0$.  The $\cu_0$ component of the Maurer-Cartan equation $d\a+\frac{1}{2}[\a, \a]=0$ gives
$$d\a_0+\frac{1}{2}([\a_0,\a_0] + [\a_1, \a_1])=0.$$
So $f^*([\a_1, \a_1])=0$.  This implies that $f^{-1}df$ is $\cu_1$-valued and the subspace $f^{-1}\Im(df_p)$ is abelian for all $p\in \co$.  This means that $(f^{-1}f_{y_1}, \ldots, f^{-1}f_{y_r})$ is a solution of the curved flat system \eqref{bd} with respect to any coordinate system $y$.  Use the Cartan-K\"ahler Theorem we can see that the curved flat system should only depend on $n$ functions of one variable,  where $n= \dim(\cu_1)-r$.
But the curved flat system \eqref{bd} is a system of $r(r-1)/2$ equations of $nr$ functions. This indicates that the curved flat system has many redundant functions and we probably can use geometry to find a special coordinate system on integral submanifolds so that their PDE involves only $n$ functions.  This is indeed the case.  We can find a special coordinate system $x$ on an integral submanifold of $(U,\ci_\sigma)$ so that the corresponding PDE  written in $x$ coordinate is gauge equivalent to the $U/U_0$-system.  

Since the curved flat system is gauge equivalent to the $U/U_0$-system,  we can use techniques from soliton theory to construct infinitely many explicit integral submanifolds and conservations laws for the exterior differential system $\ci_\sigma$.  

This paper is organized as follows: We explain the gauge equivalence of the $U/U_0$-system and the curved flat system in section 2, give a brief review of theory of exterior differential systems in section 3, and give Bryant's proof that the exterior system $\ci_\sigma$ on the Lie group $U$ is involutive in section 4.  Finally in section 5, we explain how to use the Birkhoff loop group factorization to construct infinitely many families of explicit solutions and commuting flows for the $U/U_0$-system and conservation laws for $\ci_\sigma$.  

\bs
\section{The $U/U_0$-system}

Let $G, \tau, \sigma, U, U_0$ be as in section 1,  $U/U_0$ the corresponding symmetric space, $\cu=\cu_0+\cu_1$ its Cartan decomposition, and $(\ , )$ be an ad-invariant, non-degenerate bilinear form on $\cu$. 
 
  Note that the $U/U_0$-system \eqref{ai} can be also defined invariantly as a system for maps $v:\ca\to \cu_1\cap\ca^\perp$ so that $[da, v]$ is flat, where $da$ means the differential of the identity map $a(\xi)=\xi$ on $\ca$.  When we choose a basis $a_1, \ldots, a_r$ of $\ca$, the system becomes  \eqref{ai}.  Change basis of $\ca$ amounts to linear change of coordinates of $\R^r$.  

 The $U/U_0$-system \eqref{ai} and the curved flat system \eqref{bd} are gauge equivalent.  To explain this we first recall some known Propositions, which can be proved by direct computations. 

\bprop\label{ax}
The following statements are equivalent for smooth maps $u_i:\R^n\to \cg$, $1\leq i\leq n$:
\ben
\item $\sum_{i=1}^n u_i dx_i$ is a flat $\cg$-connection $1$-form on $\R^n$,
\item the first order system 
$E_{x_i} = Eu_i, \quad 1\leq i\leq n$
is solvable,
\item there exists $g:\co\to G$ such that 
$g^{-1}dg= \sum_{i=1}^n u_i dx_i$
for some open subset $\co$ of the origin in $\R^n$.   
\een
\eprop

\bprop\label{be} (\cite{Ter97})
The following statements are equivalent for $v:\R^n\to \cu_1\cap\ca^\perp$:
\ben
\item $v$ is a solution of the $U/U_0$-system \eqref{ai},
\item $\sum_{i=1}^r [a_i, v] dx_i$ is a $\cu$-valued  flat connection on $\R^r$,
\item
\begin{equation}\label{ba}
\o_\l=\sum_{i=1}^r (a_i\l + [a_i, v])dx_i
\end{equation}
 is a $\cg$-valued  flat connection on $\R^r$ for all parameter $\l\in \C$,
\item there is an $s\in \R$ so that $\o_s= \sum_{i=1}^r (a_i s+ [a_i, v])dx_i$ is a $\cu$-valued  flat connection on $\R^r$.
\een
\eprop

\bprop \label{bf} A smooth map $(A_1, \cdots, A_r):\R^r\to \prod_{i=1}^r \cu_1$ is a solution of the curved flat system  \eqref{bd} associated to $U/U_0$  if and only if 
$$\w_\l= \sum_{i=1}^r \l A_i dx_i$$
is a flat $\cg$-valued connection $1$-form on $\R^r$ for all $\l\in \C$. 
\eprop

The flat connections $\o_\l$ and $\w_\l$  are called {\it Lax connections\/} of the $U/U_0$-system and the $U/U_0$ curved flat system.  

A map $\xi:\C\to \cg$ is said to satisfy the {\it $U/U_0$-reality condition\/} if 
$$\tau(\xi(\bar \l))=\xi(\l), \quad \sigma(\xi(\l))= 
\xi(-\l).$$
It follows from the definition that $\xi(\l)= \sum_j \xi_j\l^j$ satisfies the $U/U_0$-reality condition if and only if $\xi_j\in \cu_0$ if $j$ is even and $\xi_j\in \cu_1$ if $j$ is odd.      
Note that both Lax connections $\o_\l$ and $\w_\l$ satisfy the $U/U_0$-reality condition.

It follows from Proposition \ref{ax} that if $v$ is a solution of the $U/U_0$-system then there exists a unique $E(x,\l)$ so that 
$$\bca E^{-1}E_{x_i} = a_i\l + [a_i, v], & 1\leq i\leq r,\\
E(0,\l)=e.\eca$$
Since $\o_\l$ satisfies the $U/U_0$-reality condition, $E$ also satisfies the $U/U_0$-reality condition:
$$\tau((E(x,\bar\l))= E(x,\l), \quad \sigma(E(x,\l))= E(x, -\l).$$
We call such $E$ the {\it parallel frame\/} of the Lax connection $\o_\l$ associated to $v$.  

The following Proposition says that solutions of the $U/U_0$-system give rise to solutions of the curved flat system.

\bprop\label{bg} (\cite{Ter03})    Let $v:\R^r\to \cu_1\cap \ca^\perp$ be a solution of the $U/U_0$-system \eqref{ai},
and $E(x,\l)$ the parallel frame of the corresponding Lax connection $\o_\l$ defined by \eqref{ba}.  Let $g(x)= E(x,0)$, and $A_i= ga_ig^{-1}$ for $1\leq i\leq r$.  Then
\ben
\item[(i)] the gauge transformation of $\o_\l$ by $g$ is 
$$g\ast \o_\l= \sum_{i=1}^r \l ga_i g^{-1} dx_i,$$
\item[(ii)] $(A_1, \ldots, A_r)$ is a solution of the curved flat system \eqref{bd},
\een
\eprop

\bthm \label{bl} (\cite{FerPed96b})
  If $(A_1, \ldots, A_r)$ is a solution of the curved flat system \eqref{bd} associated to $U/U_0$, then there exists $f:\co\to U$ such that $f^{-1}f_{x_i}=A_i$ for all $1\leq i\leq r$ and  $\pi(f)$ is a curved flat in $U/U_0$, where $\pi:U\to U/U_0$ is the natural projection.  Conversely, every curved flat in $U/U_0$ can be lifted to a map $f$ to $U$ so that $(f^{-1}f_{x_i}, \ldots, f^{-1}f_{x_r})$ is a solution of the curved flat system \eqref{bd}. 
\ethm

A direct computation implies

\bprop\label{bh} 
If $(A_1, \cdots, A_r)$ is a solution of the curved flat system \eqref{bd} associated to $U/U_0$, then there exists a smooth map $f:\R^r\to U$ such that $f$ satisfies the following conditions:
\begin{equation}\label{bi}
\bca f^{-1}f_{x_i}\in \cu_1, &\\
[f^{-1}f_{x_i}, f^{-1}f_{x_j}]=0, \quad {\rm for \ all \ } i\not=j.
\eca
\end{equation}
Conversely, if $f:\R^r\to U$ is an immersion satisfying \eqref{bi}, then $$(f^{-1}f_{x_i}, \ldots, f^{-1}f_{x_i})$$ is a solution of the curved flat system \eqref{bd}.
\eprop

An immersed submanifold $f:\co\to \cu_1$ is called  {\it flat abelian \/} (\cite{Ter03})  if 
\ben
\item $[f_{y_i}, f_{y_j}]=0$ for all $ 1\leq i\not= j\leq n$,
\item the induced metric on $\co$ is flat.
\een
The following Theorems give explicit algorithms to construct flat abelian submanifolds in $\cu_1$ and curved flats in the symmetric space $U/U_0$ from solutions of the $U/U_0$-system.  The proofs can be found in \cite{Ter03}.  

\bthm \label{bb} (\cite{Ter03}) 
  Let $v$ and $E$ be as in Proposition \ref{bg}.  Set $Y=\frac{\p E}{\p \l}E^{-1} \big|\ _{\l=0}$. 
Then  $Y$ is an immersed flat abelian submanifold in $\cu_1$.  Conversely, locally all flat abelian submanifolds in $\cu_1$ can be constructed
this way.  
\ethm

  For $g\in U$ and $x\in U$, 
$$g\ast x=
gx\sigma(g)^{-1}$$
defines an action of $U$ on $U$ (it is called the $\sigma$-action).   The orbit at
$e$ is 
$$M=\{g\sigma(g)^{-1}\n g\in U\}\ \subset U.$$ 
Since the isotropy subgroup at $e$ is $U_0$, the orbit $M$ is diffeomorphic to
$U/U_0$. In fact, $M$ is a totally geodesic submanifold of $U$ and is isometric to the symmetric space $U/U_0$.  This is the classical {\it Cartan
embedding\/} of the symmetric space $U/U_0$ in $U$.

\bthm \label{bc} (\cite{Ter03})
With the same assumption as in Theorem \ref{bg}, and set
$$\psi(x)= E(x,1)E(x,-1)^{-1}.$$
Then $\psi$ is a curved flat in the symmetric space $U/U_0$.  Conversely, locally all curved flats in $U/U_0$ can be constructed this way.
\ethm

\bthm\label{bk}  Let $\co$ be an open neighborhood of $0\in\R^r$.  
If $f:\co\to U$ is an immersion satisfying \eqref{bi}, then there exists a local coordinate system $x$ near $0$, a regular basis $\{a_1, \ldots, a_r\}$ of the maximal abelian subspace $\ca= \Im(df_0)$, $g:\co\to U_0$, and   a solution $v:\R^r\to \cu_1\cap\ca^\perp$ of the $U/U_0$-system \eqref{ai} so that
\begin{equation}\label{bj}
\bca
f^{-1}f_{x_i}= ga_ig^{-1}, &\\
g^{-1}g_{x_i}= [a_i, v].
\eca
\end{equation}
Conversely, if $v$ is a solution of the $U/U_0$-system, then $f(x)= E(x,1)E(x,0)^{-1}$ satisfies \eqref{bi} and \eqref{bj}, where $E(x,\l)$ is the parallel frame for the Lax connection $\o_\l$ corresponding to $v$.  
\ethm

\begin{proof} 
Since generically all maximal abelian subalgebra are conjugate under elements of $U_0$, there exist $g:\co\to U_0$ and $b_1, \cdots, b_n:\co\to \ca$ such that 
$$ f^{-1}f_{y_i} = g b_i g^{-1}$$ 
for $1\leq i\leq n $.Ê A direct computation implies 
\begin{align*} 
0 &= d ( d f) =  d \left( f \sum_{i=1}^n gb_ig^{-1} d y_{i} \right) \\
&= f \sum_{j \neq i} \left(g(b_i)_{y_j} g^{-1} + [ g_{y_j}g^{-1}, gb_ig^{-1} ] \right) d y_j\wedge  d y_i \\ 
&= f g \left(\sum_{j \neq i} ((b_i)_{y_j} - [b_i, g^{-1}g_{y_j}]) d y_j\wedge  d y_i\right) g^{-1}. 
\end{align*} 
This implies that 
\begin{equation} \label{ag} 
(b_i)_{y_j} - [b_i, g^{-1}g_{y_j}]= (b_j)_{y_i} - [b_j, g^{-1}g_{y_i}] 
\end{equation} 
for all $i\not=j$. 
Let $(\ , )$ denote a non-degenerate ad-invariant bilinear form on $\cu$.Ê Then $(\ca,Ê~ [\ca, \cu])=0$ and $[\cu_{\ca}, [\ca, \cu])=0$.  So $[\ca, \cu]\subset \cu_{\ca}^\perp$ and $[\ca, \cu]\subset \ca^\perp$.  Also we have
$$\cu_1= \ca \oplus (\ca^\perp\cap \cu_1).$$ 
Note $(b_i)_{y_j}\in \ca$ and $[b_i, g^{-1}g_{y_j}]\in [\ca, \cu_0]$ is contained in $\ca^\perp$. By \eqref{ag} we get 
\begin{subequations}\label{ah}
\begin{gather} 
(b_i)_{y_j}=(b_j)_{y_i},\label{ah1}\\ 
[b_i, g^{-1}g_{y_i}] = [ b_j, g^{-1}g_{y_i}],\label{ah2} 
\end{gather} 
\end{subequations}
for all $1\leq i\not= j\leq n $. 
Equation \eqref{ah1} implies that $\sum_{i=1}^n b_i  d y_i$ is closed. So there exist a local coordinate change $x=x(y)$ and constant $a_1, \ldots, a_n$ in $\ca$ such that $\sum_{i=1}^n b_i  d y_i = \sum_{i=1}^n a_i  d x_i$. 

Let $\b= \sum_{i=1}^n b_i d y_i$.Ê Equation \eqref{ah2} can be rewritten asÊ $[\b, g^{-1} d g]=0$. Write $\b$ and $g^{-1} dg$ inÊ $x$ coordinate to getÊ $\b= \sum_{i=1}^n a_i dx_i$ and $g^{-1} d g= \sum_{i=1}^n g^{-1} g_{x_i}  d x_i$.Ê Then 
$$0=[\b, g^{-1} d g]= \sum_{i\not=j} [a_i, g^{-1} g_{x_j}]  dd x_i\wedge  d x_j.$$ 
So we have $$[a_i, g^{-1}g_{x_j}]= [a_j, g^{-1}g_{x_i}], \quad \forall \ i\not=j.$$ 
Ê Up to a linear change of coordinates of $x$, we may assume that $a_i$'s are regular.Ê Note the kernel of $\ad(a_i)$ on $\cu_1$ is $\ca$, and the tangent plane of the orbit $\Ad(U_0)(a_i)$ at $a_i$ is $[a_i, \cu_0]$.Ê By assumption $\Ad(U_0)(\ca)$ is open in $\cu_1$. So the  dimension of the tangent plane of the principal $\Ad(U_0)$-orbit at $a_i$  is equal to $\dim(\cu_1)-\dim(\ca)$.Ê Thus $\ad(a_i)$ maps $\ca^\perp\cap\cu_1$ isomorphically onto $\cu_0\cap (\cu_0)_{\ca}^\perp$, where $(\cu_0)_{\ca}=\{\xi\in \cu_0\n [\xi, \ca] = 0\}$.Ê Then by \eqref{ah2}, there exists a $v:\co\to \cu_1\cap\ca^\perp$ so that 
$$g^{-1}g_{x_i}= [a_i, v], \quad 1\leq i\leq n .$$ 
But $g^{-1} d g=\sum_i [a_i, v]  d x_i$ is a flat connection.  By Proposition \ref{be}, $v$ is a solution of the $U/U_0$-system \eqref{ai}. 

To prove the converse, note that 
$$E^{-1}dE=\o_\l = \sum_{i=1}^r (a_i\l+ [a_i,v]) dx_i.$$
Set $g(x)= E(x,0)$ and $F(x,\l)= E(x,\l)E(x,0)^{-1}= E(x,\l)g(x)^{-1}$.  A direct computation implies that
$$F^{-1}dF= g\o_\l g^{-1}- dg g^{-1}= \sum_{i=1}^r \l ga_i g^{-1} dx_i.$$
Note $f(x)= F(x,1)$.  So $f^{-1}df=\sum_{i=1}^r ga_i g^{-1}dx_i$.  
\end{proof} 

\brem 
The maps $g$ and $v$ in Theorem \ref{bk} are essentially unique.  To see this, suppose we have $g, \ti g$ so that 
$$ga_ig^{-1} = \ti g a_i \ti g^{-1}= f^{-1}f_{x_i},$$
 $g^{-1}g_{x_i}=[a_i, v]$, and $\ti g^{-1}\ti g_{x_i}=[a_i, \ti v]$.  Since $g^{-1}\ti g a_i=a_i$,  there exists $(\cu_0)_{\ca}$-valued map $h$ such that $g^{-1}\ti g= h$, i.e., $\ti g= gh$.  But
\begin{align*}
\ti g^{-1}\ti g_{x_i}&= [a_i, \ti v] = h^{-1}[a_i, v] h + h^{-1}h_{x_i}\\
&=[a_i, h^{-1}vh] + h^{-1} h_{x_i}\in \cu_{\ca}^\perp + \cu_{\ca}.
\end{align*}
Thus 
\begin{equation*}
\bca h^{-1} h_{x_i}=0, \\
[a_i, h^{-1}v h]= [a_i, \ti v].\eca
\end{equation*}
 The first equation implies $h$ is a constant.  Since $h^{-1}vh\in \cu_1\cap\ca^\perp$ and $\ad(a_i)$ is injective on $\cu_1\cap \ca^\perp$, the second equation implies that $h^{-1}vh = \ti v$.  This proves that $\ti g= gh$ and $\ti v= h^{-1}vh$ for some constant $h\in (\cu_0)_{\ca}$.  
\erem

\bs
\section{Basics of exterior differential systems}

We give a brief account of Cartan-K\"ahler theory based on the lectures given by R. Bryant at MSRI in 1999 and 2003 (cf. \cite{BC3G} for details and references). 

Let $M$ be a smooth manifold, and $\Omega^{*}(M)$ the graded algebra of differential forms on $M$.Ê An ideal $\ci$ of $\W^*(M)$ is called a \emph{differential ideal} if $\ci$ satisfies the following conditions: 
\begin{enumerate} 
\item $ \ci = \bigoplus_{j} \ci^{j} $, where $ \ci^{j} = \Omega^{j}(M) \cap \ci $; 
\vspace{3pt}
\item $ d \ci \subset \ci $. 
\end{enumerate} 

An \emph{exterior differential system} (EDS) is a pair $(M, \ci)$ consisting of a smooth manifold $M$ and a differential ideal $\ci \subset \Omega^{*}(M)$. 

A submanifold $N \subset M$ is called an \emph{integral submanifold} for the EDS $(M, \ci)$ if $i^{*}\ciÊ = 0$, where $i: N \hookrightarrow M $ is the inclusion. In local coordinates this condition can be written as a system of PDE (or ODE). 

A linear subspace $E \subset T_{p}M $ is said to be an \emph{integral element} of $\ci$ if $\varphi |_{E} = 0 $ for all $\varphi \in \ci$. The set of all integral elements of $\ci$ of dimension $n$ is denoted $v_{n}(\ci )$. A submanifold of $M$ is an integral submanifold of $\ci$ if and only if each of its tangent space is an integral element of $\ci$. 

Note that $v_{n}(\ci ) \cap Gr_{n}(T_{p}M)$ is a real algebraic sub-variety ofÊ $Gr_{n}(T_{p}M)$, which may be very complicated. The set of \emph{ordinary integral elements} $v_{n}^{o}(\ci) \subset v_{n}(\ci)$ consists of those which are locally cut out 'cleanly' by finite number of $n$-forms in $\ci$, so that the connected components of $v_{n}^{o}(\ci)$ are smooth embedded submanifolds of $Gr_{n}(TM)$. The rigorous definition can be found in \cite{BC3G}. 

Let $\{e_{1}, \cdots, e_{n} \}$ be a basis of the linear subspace $E$ of $T_p M$. The \emph{polar space} of $E$ is defined to be the vector space 
\[ 
H(E) = \{ v \in T_{p}M \n \varphi (v, e_{1}, \cdots, e_{n}) = 0 \textrm{ for all } \varphi \in \ci^{n+1} \} . 
\] 
When $E \in v_{n}( \ci ) $, a $(n+1)$-plane $E^{+}$ containing $E$ is an integral element of $\ci$ if and only if $E^{+} \subset H(E)$. Define 
\[ 
r(E) = \dim H(E) - \dim E -1. 
\] 
This integer may jump up at certain points. An ordinary integral element $E$ is called \emph{regular} if $r$ is locally constant in a neighborhood of $E$ in $v_{n}^{o} (\ci)$. The set of regular integral elements is denoted $v_{n}^{r}(\ci)$ and is a dense open subset of $v_{n}^{o} (\ci)$. Thus $v_{n}^{r}(\ci) \subset v_{n}^{o} (\ci) \subset v_{n}(\ci) \subset Gr_{n}(TM) $.Ê An integral submanifold is calledÊ \emph{regular} if all of its tangent spaces are regular integral elements. 

We state the following two theorems that are given in \cite{Bry99}:

\bthm \label{br} (Cartan-K\"ahler Theorem) 
Suppose $(M, \ci)$ is a real analytic EDS and that $N \subset M $ is a connected real analytic regular $n$-dimensional integral submanifold of $\ci$ with $ r(N) \geq 0 $. Let $R \subset M$ be a real analytic submanifold of codimension $r(N)$ containing $N$ such that 
\[ 
\dim (~ T_{p}R \cap H(T_{p}N) ~ ) = n+1, \textrm{ for all } p \in N . 
\] 
Then there exists a unique connected real analytic $(n+1)$-dimensional integral submanifold $\tilde{N}$ such that $N \subset \tilde{N} \subset R$ . 
\ethm 

%We will call a decomposable $n$-form $\Theta = \theta_{1} \wedge \cdots \wedge \theta_{n}$ a system of $n$ independent variables. 

A \emph{regular flag} is a flag of integral elements 
\[
 (0) = E_{0} \subset E_{1} \subset \cdots \subset E_{n} = E \subset T_{p}M 
\]
where $E_{j} \in v_{j}^{r}(\ci)$ for $0 \leq j < n$ and $E_{n} \in v_{n}(\ci)$. Note that $E_{n}$ may not be regular, but one can show that it must be ordinary. By applying Cartan-K\"ahler Theorem repeatedly to this flag, one can show that there is a real analytic $n$-dimensional integral manifold $N \subset M$ passing through $p$ and satisfying $T_{p}N = E$. Set 
\[
c(E_{j}) = \dim (T_{p}M) - \dim H(E_{j}) .
\]
\bthm (Cartan's Test)
Let $(M, \ci)$ be an EDS, and $F = (E_{0}, \cdots, E_{n}) $ an integral flag of $\ci$. Then $v_{n}(\ci)$ has codimension at least 
\[
c(F) = c(E_{0}) + \cdots + c(E_{n-1}) 
\]
in $Gr_{n}(TM) $ at $E_{n}$. Moreover, $F$ is a regular flag if and only if $v_{n}(\ci)$ is a smooth submanifold  of $Gr_{n}(TM) $ in a neighborhood $E_{n}$ and has codimension exactly $c(F)$. 
\ethm

The \emph{Cartan characters} of the flag $F$ are the numbers 
\[
s_{j}(F) = \dim H(E_{j-1}) - \dim H(E_{j}), \qquad 0 \leq j \leq n, 
\]
with the convention $c(E_{-1}) = 0 $ or $H(E_{-1}) = T_{p}M$. These numbers exhibit the generality of integral submanifolds. Roughly speaking, the integral manifolds near $N$ will depend on $s_{0}$ constants, $s_{1}$ functions of one variable, $\cdots$, $s_{n}$ functions of $n$ variables. 

A connected open subset $Z$ of $v_{n}^{o}(\ci)$ is called \emph{involutive} if every $E \in Z$ is the terminus of a regular flag. When $Z$ is clear from the context, we simply say that our EDS $(M,\ci)$ is involutive. 

Suppose $(M,\ci)$ is an EDS with $n$-dimensional integral submanifold.  A \emph{conservation law} forÊ $(M, \ci)$ is an $(n-1)$-form $\phi \in \Omega^{n-1}(M)$ such that $d (f^{*} \phi) = 0 $ for every integral submanifold $f: N^{n} \hookrightarrow M$ of $\ci$. Actually, one only considers as conservation laws those $\phi$ such that $d \phi \in \ci$. Two ``trivial'' type of conservation laws are $\phi \in \ci^{n-1}$ or $\phi$ being exact on $M$. Factoring outÊ these cases, \emph{the space of conservation laws} is defined to be $\mathcal{C} = H^{n-1}(\Omega^{*}(M) / \ci )$. It also makes sense factor out those conservation laws represented by closed $(n-1)$-forms on $M$ (then the quotient space is called the space of \emph{proper} conservation laws). One can study the symmetries of the EDS and then apply Noether's Theorem to compute the corresponding conservation laws (cf. \cite{BCG} for details).

\bs
\section{Involutivity of the EDS $(U, \ci_\sigma)$} 

Let $G, \tau, \sigma, U, \cu_0$ and $\cu_1$ be as in section 1. 
Let $\alpha$ be the canonical left invariant $1$-form $g^{-1}dg$ on $U$. Write 
$$\alpha = \alpha_{0} + \a_1$$ 
with respect to the Cartan decomposition $\cu=\cu_0+\cu_1$. The $\cu_j$-componentÊ of the Maurer-Cartan equation $d \alpha + \frac{1}{2} [\alpha, \alpha] = 0$ gives 
$$ \bca d\a_0+ \frac{1}{2}([\a_0, \a_0]+[\a_1,\a_1])=0,&\\
d\a_1+ [\a_0, \a_1] =0.
\eca$$ 

Let $\ci_\sigma \subset \Omega^{*}(U) $ be the differential ideal generated by the components of $\alpha_{0}$. It follows from the Maurer-Cartan equation that 
\begin{eqnarray*} 
\ci_\sigma & = & \langle \alpha_{0}, d \alpha_{0}\rangleÊ \\ 
& = & \langle \alpha_{0},  [ \alpha_{1}, \alpha_{1} ] \rangle . 
\end{eqnarray*} 
Here $\langle \ ,  \rangle $ denotes the algebraic ideal generated by the enclosed forms. 

The following Proposition was proved by R. Bryant. 
\bprop \label{ac} (\cite{Bry03})
The EDS $(U, \ci_\sigma)$ is involutive. 
\eprop 

\begin{proof} 
Since everything is homogeneous, we only need to look at the integral elements $E \subset T_{e}U = \cu $. Note that $E \subset \cu_{1} = \cap_{j \neq 1} \ker (\alpha_{j}) $. For $E=(0) \in v_{0}(\ci_\sigma)$, we have $H(E) = \cu_{1}$ and $v_{0}(\ci_\sigma) \cong U$. Thus $v_{0}(\ci_\sigma) = v_{0}^{o}(\ci_\sigma) = v_{0}^{r}(\ci_\sigma) $. Now consider $E = \R x \in v_{1}(\ci_\sigma)$ for some $x \in \cu_{1} - \{ 0 \}$. Its polar space is 
$$ H(E) = \{ y \in \cu_{1} | [x, y] = 0 \} ,$$ 
since $[ \alpha_{1}, \alpha_{1} ]_{e}(x,y) = [x,y] $. For generic such $x$, $H(\R x)$ will be a maximal abelian subalgebra of $\cu_{1}$, and set $\dim H(\R x) = \dim \ca = r $. Therefore 
$$ v_{1}(\ci_\sigma) = v_{1}^{o}(\ci_\sigma) \varsupsetneqq v_{1}^{r}(\ci_\sigma) .$$ 
Furthermore, when $\R x \in v_{1}^{r}(\ci_\sigma)$, every subspace $E$ of $H(\R x)$ containing $\R x$ is also regular and has $ H(E) = H(\R x ) $. Thus generic $E \in v_{r}^{o}(\ci_{\sigma})$ is the terminus of a regular flag, and our EDS is involutive. In fact, every regular integral curve of $\ci_\sigma$ lies in a unique $r$-dimensional integral submanifold of $\ci_\sigma$,  or locally the integral submanifolds depend on $s_{0}= \dim \cu - \dim \cu_{1}$ constants and $s_{1} = \dim (\cu_{1}) - r $ functions of one variable (since $s_{2}=\cdots = s_{r}=0$). 
\end{proof}

\bcor\label{ae}  Let  $\co$ be an open subset of $\R^r$, and  $f:\co\to U$ an immersion.  Then the following statements are equivalent:
\ben
\item $f$ is a $k$-dimensional integral submanifold of $(U,\ci_\sigma)$,
\item $f$ satisfies \eqref{bi},
\item $(f^{-1}f_{x_1}, \ldots, f^{-1}f_{x_r})$ is a solution of the curved flat system \eqref{bd}.
\een
\ecor

Hence the PDE for the EDS $(U, \ci_\sigma)$ is the curved flat system \eqref{bd} associated to $U/U_0$.  

As a consequence of the Cartan K\"ahler Theorem 
\ref{br}, Proposition \ref{ac} and Corollary \ref{ae} that real analytic curved flats in $U/U_0$ or real analytic solutions of the curved flat systems \eqref{bd} depend only on $\dim(\cu_1\cap \ca^\perp)$ functions of one variable along a non-characteristic line.  

By Theorem \ref{bk}, there is a special coordinate system $x$ on $\R^r$ so that the curved flat system \eqref{bd} written in $x$ coordinate system is gauge equivalent to the $U/U_0$-system \eqref{ai}. The Cartan-K\"ahler theory implies that the Cauchy problem of the $U/U_0$-system has a unique local solution for any given local real analytic initial data on the $x_1$-axis.  But it was also proved in \cite{Ter97} using the inverse scattering theory of Beals and Coifman (\cite{BeaCoi84}) that given any smooth rapidly decaying function $v_0:\R\to \cu_1\cap\ca^\perp$, there exists a unique smooth solution $v:\R^r\to \cu_1\cap\ca^\perp$ so that 
$$v(x_1, \ldots, x_r)= v_0(x_1, 0,\ldots, 0).$$ 
 Although the theory exterior differential system seems to give a weaker result concerning the Cauchy problem, it may prove to be a very good tool to detect ``integrablility''. 

\brem  Let $G, \tau, U$ be as above, and $\rho$ an order $k$ automorphism of $G$ so that $d\rho_e$ is complex linear.  Assume that 
$$\tau\rho= \rho^{-1}\tau^{-1}.$$
Let $\cg_j$ denote the eigenspace of $d\rho_e$ with eigenvalue $e^{\frac{2\pi {\rm i\/} j}{k}}$, and $\ck_j= \cu\cap \cg_j$.  Then 
\begin{equation}\label{ab}
\cu= \ck_0+ \cdots + \ck_{k-1}.
\end{equation}
Let $\a=g^{-1}dg$, and 
$$\a=\a_0+\cdots + \a_{k-1}$$
the decomposition of $\a$ with respect to \eqref{ab}.  Let $\ci_\rho$ denote the differential ideal on $U$ generated by $\a_0, \a_2, \ldots, \a_{k-1}$.  Then 
\begin{align*}
\ci_\rho &=\li \a_0, \a_2, \ldots, \a_{k-1}, d\a_0, d\a_2, \ldots, d\a_{k-1}\ri \\
&= \li \a_0, \a_2, \ldots, \a_{k-1}, [\a_1, \a_1]\ri. 
\end{align*}
We define regular elements in $\ck_1$ the same way as before, namely, $a\in \ck_1$ is regular if  it is contained in a maximal abelian subspace  $\ca$ in $\ck_1$ and $\Ad(U_0)(\ca)$ is open in $\ck_1$.  If $\ck_1$ admits regular elements, then the  proof of Proposition \ref{ac} works for $(U,\ci_\rho)$.  In fact, in this case, we have:
\ben
\item $(U, \ci_\rho)$ is involutive.
\item If $\dim(\ca)=r$, then any $r$-dimensional integral submanifold depend on $\dim(\ck_1)-\dim(\ca)$ number of functions of one variable.
\item every regular integral curve is contained in a unique $r$-dimensional integral submanifold of $(U,\ci_\rho)$. 
\item The curved flat system associated to $U/K$ is the system \eqref{bd} for $(A_1, \ldots, A_r):\R^r\to \ck_1$, and the $U/K$-system is the system \eqref{ai} for $v:\R^r\to \ck_{k-1}$. Modulo a change of coordinate system of $\R^r$, these two system are gauge equivalent.  
\item Given an immersion $f:\R^r\to U$, the following statements are equivalent:
\ben
\item $f$ is an integral submanifold of $(U,\ci_\rho)$, 
\item $f^{-1}f_{x_i}\in \ck_1$ and $[f^{-1}f_{x_i}, f^{-1}f_{x_j}]=0$ for all $i, j$,
\item $(f^{-1}f_{x_1}, \ldots, f^{-1}f_{x_r})$ is a solution of the curved flat system associated to $U/K$. 
\een
\item The PDE for the EDS $(U,\ci_\rho)$ is the curved flat system associated to $U/K$.
\een
\erem

\bs
\section{Conservation laws and commuting flows\/}

We construct infinitely many conservation laws and commuting flows for the $U/U_0$-system, and indicate how to construct infinitely many explicit solutions of the $U/U_0$-system.  

First we review the Birkhoff Factorization Theorem (for detail see \cite{PreSeg86}).  Let $\e>0$ be a small number, and $\co_\e=\{\l\in \C\n \frac{1}{\e}<|\l |\leq \infty\}$ the open neighborhood at $\infty$ in $S^2=\C\cup\{\infty\}$.  $L(G)$ denote the group of holomorphic maps $g:\co_\e\setminus\{\infty\} \to G$, $L_+(G)$ the subgroup of $g\in L(G)$ such that $g$ can be extended to a holomorphic map  in $\C$, and $L_-(G)$ the subgroup of $g\in L(G)$ that can be extended to a holomorphic map in $\co_\e$ and is equal to the identity $e$ at $\infty$.  

\bthm (Birkhoff Factorization Theorem)  The multiplication map 
$$\mu: L+(G)\times L_-(G)\to L(G), \quad (g_+, g_-)\mapsto g_+g_-$$
is one to one, and the image is an open dense subset of $L(G)$. 
\ethm

In other words, for generic $g\in L(G)$, we can factor $g= g_+g_-$ uniquely with $g_\pm \in L_\pm(G)$. Let $\hat e$ denote the constant map from $\co_e\setminus \{\infty\}$ to $G$ with constant $e$. Since $\hat e$ lies in the image of the multiplication map $\mu$, there is an open subset of $\hat e$ so that all elements in this open subset can be factored uniquely.  

Let $\ti\tau$ and $\hat \sigma$ denote the map on $L(G)$ defined by
$$(\ti \tau(g))(\l)= \tau(g(\bar \l)), \quad (\hat \sigma(g))(\l)= \sigma(g(-\l)).$$
It is easy to check that 
\ben
\item
$\hat\tau$ and $\hat\sigma$ are conjugate linear and complex linear involutions of $L(G)$,
\item $g\in L(G)$ is a fixed point of both $\ti\tau$ and $\hat\sigma$ if and only if $g$ satisfies the $U/U_0$-reality condition: $\tau(g(\bar\l))= g(\l)$, $\sigma(g(\l))=\sigma(-\l)$. 
\item  both $\ti\tau$ and $\hat\sigma$ leave $L_\pm(G)$ invariant. 
\een
  Let $L^{\tau,\sigma}(G)$ and $L_\pm^{\tau,\sigma}(G)$ denote the subgroup of fixed points of $\ti\tau$ and $\hat\sigma$ of $L(G)$ and $L_\pm(G)$ respectively.  Then we have

\bcor \label{aj}
The multiplication map 
 $$L_+^{\tau,\sigma}(G)\times L_-^{\tau,\sigma}(G)\to L^{\tau,\sigma}(G)$$
 is one to one and the image is open and dense in $L^{\tau,\sigma}(G)$.   \ecor

We want to use this factorization to construct infinitely many solutions and commuting flows for the $U/U_0$-system.  Given $b\in \ca$ and $j>1$ an odd integer, $x\in \R^r$, and $t\in \R$, let $e^A(x,t)\in L_+^{\tau,\sigma}(G)$ be defined by 
$$e^{A}(x,t)(\l)= \exp((a_1x_1+\cdots +a_rx_r)\l + b\l^j t).$$
Given $f\in L_-^{\tau,\sigma}(G)$, since $e^A(0,0)=\hat e$ is the identity in $L^{\tau,\sigma}(G)$ and $e^A$ is smooth from $\R^r\times \R$ to $L^{\tau,\sigma}(G)$, by Corollary \ref{aj} there is an open subset of $(0,0)$ in $\R^r\times \R$ so that 
we can factor $f^{-1}e^A(x,t)$ uniquely as
\begin{equation}\label{ak}
f^{-1}e^A(x,t)= E(x,t)m(x,t)^{-1},
\end{equation}
where $E(x,t)\in L_+^{\tau,\sigma}(G)$ and $m(x,t)\in L_-^{\tau,\sigma}(G)$.  

Given $c\in \ca$, let 
\begin{equation}\label{al}
m^{-1}cm = Q_{c,0} + Q_{c,1} \l^{-1} +Q_{c,2}\l^{-2} +\cdots 
\end{equation}
denote the Taylor series of $(m(x,t)^{-1}cm(x,t))(\l)$ at $\l=\infty$. 
Since $m(x,t)(\l)$ $ = e$ at $\l=\infty$, 
\begin{equation}\label{au}
Q_{c,0}=c.
\end{equation}  

We want to explain how to compute $Q_{c,n}$, $m^{-1}dm$ and $E^{-1}dE$ next.  To do this, we take $\p_{x_i}$ of \eqref{ak} to get
$$f^{-1}e^A(x,t) a_i\l = E_{x_i} m^{-1} - Em^{-1}m_{x_i}m^{-1}.$$
Multiply $E^{-1}$ on the left and $m$ on the right of the above equation and use \eqref{ak} to  get
\begin{equation}\label{am}
m^{-1}a_i\l m= E^{-1}E_{x_i} - m^{-1}m_{x_i}.
\end{equation}
Take $\p_t$ of \eqref{ak} and use a similar calculation to get
\begin{equation}\label{an}
m^{-1}b\l^j m= E^{-1}E_{t} - m^{-1}m_{t}.
\end{equation}
Note that $E^{-1}E_{x_i}$ and $E^{-1}E_t$ lie in the Lie algebra $\cl_+^{\tau,\sigma}(\cg)$ of $L_+^{\tau,\sigma}(G)$, and $m^{-1}m_{x_i}$ and $m^{-1}m_t$ lie in the Lie algebra $\cl_-^{\tau,\sigma}(\cg)$ of $L_-^{\tau,\sigma}(G)$.  But it follows from the factorization theorem that
$$\cl(\cg) =\cl_+^{\tau,\sigma}(\cg) \oplus \cl_-^{\tau,\sigma}(\cg)$$ as direct sum of vector spaces.  Let $\xi_\pm$ denote the $\cl_\pm^{\tau,\sigma}(\cg)$ component of $\xi\in \cl^{\tau,\sigma}(\cg)$.    Then \eqref{am} and \eqref{an} imply that
\begin{subequations}\label{as}
\begin{gather}
E^{-1}E_{x_i}= (m^{-1}a_im\l)_+, \label{as1}\\
E^{-1}E_t =(m^{-1}bm\l^j)_+ \label{as2}\\
m^{-1}m_{x_i}= -(m^{-1}a_im\l)_-, \label{as3}\\
m^{-1}m_t= -(m^{-1}bm\l^j)_-. \label{as4}
\end{gather}
\end{subequations}

 Use \eqref{al} to see that
\begin{align*}
(m^{-1}a_im\l)_+ &= Q_{a_i, 0}\l + Q_{a_i, 1}, \\
(m^{-1}bm\l^j)_+ &= Q_{b,0}\l^j + Q_{b,1}\l^{j-1} + \cdots + Q_{b,j}.
\end{align*} 
So we get
\begin{equation}\label{ao}
\bca
E^{-1}E_{x_i}= Q_{a_i, 0}\l + Q_{a_i, 1}, &\\
E^{-1}E_t = Q_{b,0}\l^j + Q_{b,1}\l^{j-1} + \cdots + Q_{b,j}.
\eca
\end{equation}

\blem
If $c_1, c_2\in \ca$, then 
\begin{subequations}\label{ap}
\begin{gather}
[m^{-1}c_1m, \ m^{-1}c_2m]=0, \label{ap:1}\\
[m^{-1}c_1m, \ -(m^{-1}c_2\l^nm)_-]= [m^{-1}c_1m, \ (m^{-1}c_2\l^nm)_+]. \label{ap:2}
\end{gather}
\end{subequations}
\elem

\bthm (\cite{Ter97})
There exists $v:\R^r\times\R\to \cu_1\cap\ca^\perp$ so that 
$$Q_{a_i, 1}= [a_i, v].$$ Moreover, for each $t\in \R$, $v(\cdots, t)$ is a solution of the $U/U_0$-system.
\ethm

\begin{proof}
By \eqref{ap:1}, 
$$[m^{-1}a_im, \ m^{-1}a_jm]=0, \quad 1\leq i\not=j\leq r.$$
So the coefficient of $\l^{-1}$ of the left hand side has to be zero, i.e., 
$$[a_i, Q_{a_j, 1}] +[Q_{a_i,1}, a_j]=0, \quad 1\leq i\not=j\leq r.$$
But this implies that there exists $v:\R^r\times\R\to \cu_1\cap\ca^\perp$ so that 
$$Q_{a_i, 1}= [a_i, v].$$
By \eqref{ao} and Proposition \ref{ax} we see that $\sum_i (a_i\l +[a_i, v])dx_i$ is a flat $\cg$-valued connection on $\R^r$ for all $\l\in \C$.  Hence for each fixed $t$, $v(\cdots, t)$ is a solution of the $U/U_0$-system. 
\end{proof}

The following is well-known (cf. \cite{Sat84, TerUhl98}):

\bthm
\ben
\item $Q_{b,j}(x,t)$ is a polynomial in $u, \p_x v, \cdots, \p_x^{j-1}v$, 
\item $Q_{b,j}$ satisfies the following recursive formula
\begin{equation}\label{ar}
(Q_{b,j})_{x_i} + [[a_i, v], Q_{b,j}]= [Q_{b,j+1}, a_i],
\end{equation}
\item $Q_{b,0}=b$, $Q_{b,1}=[b,v]$. 
\een 
\ethm

\begin{proof}  A direct computation gives
\begin{align*}
(m^{-1}bm)_{x_i} &= [m^{-1}bm, m^{-1}m_{x_i}], \quad {\rm by \ \eqref{as1}}\\
&= [m^{-1}bm, -(m^{-1}a_i\l m)_-], \quad {\rm by\  \eqref{ap:2}}\\
&= [m^{-1}bm, (m^{-1}a_i\l m)_+].
\end{align*}
Substitute \eqref{al} to the above equation to get
\begin{equation}\label{aq}
(m^{-1}bm)_{x_i} = [m^{-1}bm, a_i\l+ u].
\end{equation}
Compare coefficient of $\l^{-j}$ of $\l^{-j}$ of \eqref{aq} to get \eqref{ar}.  

It was proved in \cite{Sat84, TerUhl98} that $Q_{b,j}$ is a polynomial in $u_i, \p_{x_i}u_i, \ldots, \p_{x_i}^{j-1}u_i$, where $u_i=[a_i,v]$.  Since $\ad(a_i)$ is a linear isomorphism between $\cu_1\cap\ca^\perp$ and $\cu_0\cap \cu_{\ca}^\perp$, (1) follows.  Since $m(\cdots, \infty)=\I$, $Q_{b,0}=b$. Use \eqref{ar} to prove  $Q_{b,1}=[b,v]$.  
\end{proof}

Use \eqref{ao} and Proposition \ref{ax} to see that 
$$\Theta^{b,j}_\l= \sum_{i=1}^r (a_i\l+ [a_i,v])dx_i + (b\l^j + Q_{b,1}\l^{j-1}+ \cdots +Q_{b,j})dt$$
is a flat connection on $\R^r\times \R$ for all $\l\in \C$.  It follows from the recursive formula \eqref{ar} and the flat equation 
$$d\Theta^{b,j}_\l +\Theta^{b,j}_\l\wedge \Theta^{b,j}_\l=0$$ 
that we have
\begin{equation}\label{ay}
\bca
[a_i,v_{x_j}]=[a_j, v_{x_i}] + [[a_i,v], [a_j,v]], \quad i\not=j, \\
[a_i, v_t] = (Q_{b, j})_x + [[a_i, v], Q_{b,j}], \quad 1\leq i\leq r.
\eca
\end{equation}
The first set of equations just means $v(\cdots, t)$ is a solution of the $U/U_0$-system for each $t$, and the second set of equations give the flow on the space of solutions of the $U/U_0$-system.   

Let $A_\C$ denote the subgroup of $G$ whose Lie subalgebra is $\ca\otimes \C$, and $L^{\tau,\sigma}_+(A_\C)$ is the subgroup of $f\in L_+^{\tau, \sigma}(G)$ such that $g(\l)\in A_\C$ for all $\l\in \C$.  Given $b\in \ca$ and $j$ a positive integer, then $\xi_{b,j}$ lies in the Lie algebra $\cl^{\tau,\sigma}_+(\ca\otimes\C)$, where  $\xi_{b,j}(\l)= b\l^j$.  Let $e_{b,j}(t)$ be the one-parameter subgroup in $L_+^{\tau,\sigma}(G)$ generated by $\xi_{b,j}$, i.e., 
 $$e_{b,j}(t)(\l)= e^{b\l^j t}.$$
It was proved in \cite{TerUhl98} that if $v(x,t)$ is a solution of \eqref{ay} then
 $$e_{b,j}(t)\cdot v(\cdots, 0):= v(\cdots, t)$$
is the dressing action of $e_{b,j}(t)\in L_+^{\tau,\sigma}(A_{\C})$  on the space of solutions of the $U/U_0$-systems.  The second set of equations of \eqref{ay} is the vector field on the space of solutions of the $U/U_0$-system corresponding to the one-parameter subgroup generated by $\xi_{b,j}$.  Since the group $L^{\tau,\sigma}_+(A_\C)$ is abelian, the flows generated by these $\xi_{b,j}$ are commuting.    
So we have 

\bthm \label{bn} (\cite{Ter97}, \cite{TerUhl98})  Given $b\in\ca$ and a positive integer $j$, the flow
\begin{equation}\label{bo}
[a_i, v_t] = (Q_{b,j})_{x_i} + [[a_i, v], Q_{b,j}], \quad 1\leq i\leq r,
\end{equation}
leaves the space of solutions of the $U/U_0$-system \eqref{ai} invariant.  Moreover, all these flows commute. 
\ethm

We sketch the method of  constructing solutions of the $U/U_0$-system below.
Let $e^{A_0}(x)\in L_+^{\tau,\sigma}(G)$ be defined by
$$e^{A_0}(x)= \exp\left({\sum_{i=1}^r a_i x_i \l}\right).$$ 

\bthm \label{} (\cite{TerUhl00a})
 Given $f\in L_-^{\tau,\sigma}(G)$, factor
\begin{equation}\label{bq}
f^{-1}e^{A_0(x)} = E(x)m^{-1}(x)
\end{equation}
with $E(x)\in L_+^{\tau,\sigma}(G)$ and $m(x)\in L_-^{\tau,\sigma}(G)$.
Expand $m(x)(\l)$ at $\l=\infty$:
$$m(x)(\l)= e+ m_{-1}(x)\l^{-1} + m_2(x)\l^{-2} + \cdots.$$
Then 
\ben
\item $m_{-1}(x)\in \cu_1$,
\item $v=(m_{-1})^\perp$ is a solution of the $U/U_0$-system, where $(m_{-1})^\perp$ is the projection of $m_{-1}$ onto $\cu_1\cap \ca^\perp$ with respect to $\cu_1 =\ca \oplus (\cu_1\cap\ca^\perp)$.
\een
\ethm 

\begin{proof}
  Use the same computation as for the proof of \eqref{as1} to conclude that 
$$E^{-1}E_{x_i}= (m^{-1}a_i m)_+ = a_i\l + Q_{a_1, 1}.$$
Expand $m(x)(\l)$ at $\l=\infty$:
$$m(x)(\l)= e+ m_{-1}(x)\l^{-1} + m_{-2}(x)\l^{-2} + \cdots.$$
A direct computation implies that 
$$m^{-1}a_i m= a_i + [a_i, m_{-1}]\l^{-1} + \cdots.$$
Therefore $Q_{a_i, 1}= [a_i, m_{-1}]$. Since $m\in L^{\tau,\sigma}(G)$, $m_{-1}(x)\in \cu_1$.  So 
$$[a_i, v]=[a_i, m_{-1}^\perp]=[a_i,m_{-1}]= Q_{a_i, 1}.$$
Hence we have shown that 
$$E^{-1}E_{x_i}= a_i\l + [a_i, v], \quad 1\leq i\leq r.$$
By Proposition \ref{be}, $v$ is a solution of the $U/U_0$-system.
\end{proof}

\brem  It was proved in \cite{TerUhl00a} that if each entry of $f\in L_-^{\tau,\sigma}(G)$ is a meromorphic function on $S^2=\C\cup\{\infty\}$, then the factorization \eqref{bq} can be carried out explicitly using residue calculus. In particular, $m(x)(\l)$ and $E(x,\l)=E(x)(\l)$ can be given by explicit formulas.  Therefore, we get explicit solutions $v=(m_{-1})^\perp$ for the $U/U_0$-system.  Since the parallel frame $E(x,\l)$ for the solution $v$ is also given explicitly, it follows from Corollary \ref{ae} and Theorem \ref{bk} that $F(x)= E(x,1)E(x,0)^{-1}$ is an explicit integral submanifold of the EDS $(U,\ci_\sigma)$.    
\erem

\ms

Next we derive conservation laws of the flows for the $U/U_0$-system.

\bthm
Let $c, b\in \ca$, and $n$ a positive integer.  Then
\begin{equation}\label{az}
(Q_{c,n}, a_i)_{x_j} = (Q_{c,n}, a_j)_{x_i}, \quad 1\leq i\not=j\leq r.
\end{equation}
In particular,  
\begin{equation}\label{bm}
\phi_{c,n}:= \sum_{i=1}^r (Q_{c,n}, a_i)\ dx_i
\end{equation}
is a closed $1$-form on $\R^r$
\ethm

\begin{proof}
Compute directly to get
\begin{align*}
(m^{-1}cm, a_i)_{x_j} &=([m^{-1}cm, m^{-1}m_{x_j}], a_i), \ \ {\rm by \ \eqref{as3}\/},\\
&= ([m^{-1}cm, -(m^{-1}a_j\l m)_-], a_i), \ \ {\rm by \ \eqref{ap:2}},\\
&= ([m^{-1}cm, (m^{-1}a_j\l m)_+], a_i) = ([m^{-1}cm, a_j\l + Q_{a_j, 1}], a_i).
\end{align*}
Use \eqref{al} and compare coefficient of $\l^{-n}$ of the above equation to get
\begin{align*}
(Q_{c,n}, a_i)_{x_j} &= ([Q_{c,n}, Q_{a_j,1}], a_i) + ([Q_{c,n+1},a_j], a_i)\\
& = (Q_{c,n}, [Q_{a_j, 1}, a_i]) +(Q_{c,n+1}, [a_j, a_i]) = (Q_{c, n}, [[a_j, v], a_i]) + 0\\
&= (Q_{c,n}, [[a_i, v], a_j]) = ([Q_{c,n}, [a_i, v]], a_j), \ \ {\rm by \ \eqref{ar}\/},\\
&=((Q_{c,n})_{x_i} -[Q_{n+1, a_i}], \ a_j) = ((Q_{c,n})_{x_i}, \ a_j).
\end{align*}
\end{proof}

\ms

If $f:\co\to U$ is a $r$-dimensional integral submanifold of the EDS $(U,\ci_\sigma)$, then by Theorem \ref{bk} and Corollary \ref{ae} there exist a special local coordinate system $x$ of $\co$, $g:\co\to U_0$ and a solution $v$ of the $U/U_0$-system \eqref{ai} such that $f^{-1}f_{x_i}= ga_i g^{-1}$ and $g^{-1}g_{x_i}=[a_i, v]$ for all $1\leq i\leq r$.  Let $\ast$ denote the Hodge star operator for the Euclidean space $\R^r$.  Given $1\leq i\not= j\leq r$, let
$$\psi_{c,n}^{ij}= \phi_{c,n}\wedge (\ast(dx_i\wedge dx_j))= \left(\sum_{\ell=1}^r (Q_{c,n}, a_\ell) dx_\ell\right)\wedge(*(dx_i \wedge dx_j)).$$
Then $\psi_{c,n}^{ij}$ is a closed $(r-1)$-form on the integrable submanifold.  In other words, $\psi^{ij}_{c,n}$ is a conservation law for the EDS $(U, \ci_\sigma)$ for all $1\leq i< j\leq r$, $c\in \ca$, and positive integer $n$.  

\ms

Next we derive the conservation laws for the flow \eqref{bo} on the space of solutions of the $U/U_0$-system \eqref{ai}.  
Given $a, c\in \ca$, compute
\begin{align*}
(m^{-1}cm, a)_t &= ([m^{-1}cm, m^{-1}m_t],a), \quad {\rm by \ \eqref{as4}\/}\\
&= ([m^{-1}cm, -(m^{-1}b\l^jm)_-], a), \quad {\rm by\ \eqref{ap:2}\/}\\
&=([m^{-1}cm, (m^{-1}b\l^j m)_+], a).
\end{align*}
Substitute \eqref{al} to the above equation and compare coefficient of $\l^{-n}$ to get
\begin{equation}\label{at}
((Q_{c,n})_t, a) = \sum_{i=0}^{j-1} ([Q_{c,n+i}, Q_{b,j-i}], a).
\end{equation}
Here we have used 
$$([Q_{c,n}, Q_{b,0}], a)= ([Q_{c,n}, b], a)= (Q_{c,n}, [b,a])=0.$$

We claim that 
\begin{equation}\label{av}
([Q_{c,n}, Q_{b,j}], a_i)= \sum_{i=1}^j (Q_{c,n+i-1}, Q_{b, j-i})_{x_i}
\end{equation}
We prove this claim by induction on $j$.  For $j=1$, we have
\begin{align*}
(Q_{c,n}, Q_{b,1}], a_i)&= (Q_{c,n}, [Q_{b,1}, a_i]) = (Q_{c,n}, [[b,v], a_i]),\\
&=(Q_{c,n}, [[a_i,v], b])= -([[a_i, v], Q_{c,n}], b)\\
& = -([Q_{c,n+1},a_i] - (Q_{c,n})_{x_i}, b)  = ((Q_{c,n})_{x_i}, b).
\end{align*}
(We used the Jacobi identity for the first line of the computation above).
This proves \eqref{av} for $j=1$.  Now assume \eqref{av} is true for $j$ and we want to prove the identity for $j+1$. We compute
\begin{align*}
([Q_{c,n}, Q_{b,j+1}], a_i) &= (Q_{c,n}, [Q_{b,j+1}, a_i]), \quad {\rm by \ \eqref{ar}\/}\\
&= (Q_{c,n}, (Q_{b,j})_{x_i} + [u_i, Q_{b,j}]) \\
&= (Q_{c,n}, Q_{b,j})_{x_i} -((Q_{c,n})_{x_i}, Q_{b,j}) + (Q_{c,n}, [u_i, Q_{b,j}])\\
&= (Q_{c,n}, Q_{b,j})_{x_i} -((Q_{c,n})_{x_i}, Q_{b,j}) -([u_i, Q_{c,n}], Q_{b,j})\\
&= (Q_{c,n}, Q_{b,j})_{x_i}  + ([a, Q_{c,n+1}], Q_{b,j}), \quad {\rm by\ \eqref{ar},\/} \\
& = (Q_{c,n}, Q_{b,j})_{x_i}  + (a, [Q_{c,n+1}, Q_{b,j}])
\end{align*}

Then the induction hypothesis implies \eqref{av} is true for $j+1$.  

It follows from \eqref{at} and \eqref{av} that we have

\bthm\label{bp}
 Let $v:\R^r\times \R\to \cu_1\cap {\ca}^\perp$ be a solution of \eqref{ay}, $c\in\ca$, and $n$ a positive integer.  Then 
\begin{equation}\label{aw}
(Q_{c,n}, a_i)_t = \sum_{\ell=0}^{j-1} \sum_{s=1}^{j-\ell}(Q_{c, n+\ell+s-1}, Q_{b, j-\ell-s})_{x_i}.
\end{equation}
\ethm

As a consequence, we  see that  $$\int_{\R^r} (Q_{c,n}, a_i) dx_1\wedge \cdots \wedge dx_r$$ is a conserved quantity for the flow \eqref{bo} on the space of rapidly decaying solutions of the $U/U_0$-system.
%%%%%%%%%%%%%%%%%%%%%%%%%%%%%%%%%%%%%%%%%%%%

\bibliographystyle{alpha}

\end{document}